# SADDLEPOINT APPROXIMATION FOR STUDENT'S $T$-STATISTIC WITH NO MOMENT CONDITIONS

By Bing-Yi Jing[1], Qi-Man Shao[2] and Wang Zhou

*Hong Kong University of Science and Technology, University of Oregon and National University of Singapore, and National University of Singapore*

A saddlepoint approximation of the Student's $t$-statistic was derived by Daniels and Young [*Biometrika* **78** (1991) 169–179] under the very stringent exponential moment condition that requires that the underlying density function go down at least as fast as a Normal density in the tails. This is a severe restriction on the approximation's applicability. In this paper we show that this strong exponential moment restriction can be completely dispensed with, that is, saddlepoint approximation of the Student's $t$-statistic remains valid without any moment condition. This confirms the folklore that the Student's $t$-statistic is robust against outliers. The saddlepoint approximation not only provides a very accurate approximation for the Student's $t$-statistic, but it also can be applied much more widely in statistical inference. As a result, saddlepoint approximations should always be used whenever possible. Some numerical work will be given to illustrate these points.

**1. Introduction.** In many statistical applications approximations to the probability that a random variable (r.v.), say $T_n$, exceeds a certain threshold value are important since the exact distribution function (d.f.) of $T_n$ may be very difficult or even impossible to obtain in most cases. Such approximations are useful, for example, in constructing confidence intervals and in calculating $p$-values in hypothesis testing. In those circumstances, we are usually dealing with tail probabilities of the r.v., $T_n$. Since these tail

Received May 2003; revised February 2004.

[1]Supported in part by Hong Kong RGC Competitive Earmarked Research Grant HKUST6148/01P.

[2]Supported in part by NSF Grant DMS-01-03487 and Grants R-146-000-038-101 and R-1555-000-035-112 at the National University of Singapore.

probabilities are typically small, accurate approximations are particularly important.

The "naive" method is to use the *Normal approximation*, which holds under mild conditions. However, this approximation is often too rough to be useful for small to moderate sample sizes. A more refined approximation is the *Edgeworth expansion* under some extra conditions. In general, the Edgeworth expansion improves the Normal approximation, but can still be inaccurate in the tails.

To overcome the difficulties encountered by the Normal approximation and the Edgeworth expansion, one can consider using a *saddlepoint approximation*, which provides a very good approximation to the tail, as well as in the center of the distribution. By a "good" approximation, here, we imply one with a small relative error. By comparison, the Edgeworth expansion gives only *absolute errors*. However, when dealing with tail probabilities, the relative error behavior is more important than the absolute error behavior. For instance, an error of 0.005 is of little importance when considering tests of size 0.05, but is of great importance when considering tests of size 0.01. Put another way, if the true probability is 0.01, it is not of much help to know that the approximation has absolute error of size $O(n^{-1})$ when $n$ is smaller than, say, 100. When, instead, the relative error is $O(n^{-1})$, we have a much more useful statement. It is quite common in statistical practice to consider test probabilities of the order of 1%, but even smaller probabilities are of interest in certain test situations. If, for example, one wishes to investigate whether a chemical substance causes cancer, one will be interested in very small test probabilities to make a convincing case. In other fields, such as reliability, small probabilities are the rule rather than the exception.

Saddlepoint approximations have been widely studied and used in many areas in recent years due to their excellent performance. For more details on the statistical importance and applications of saddlepoint approximations, one can refer to the books by Field and Ronchetti (1990), Kolassa (1997), Jensen (1995), Davison and Hinkley [(1997), Section 9.5] and to the excellent review paper by Reid (1988). All the literature clearly demonstrates how remarkably accurate the saddlepoint approximation can be. Accordingly, one should always use it if it is available.

It is worth mentioning that the extreme accuracy of the saddlepoint approximation is achieved at a cost of requiring a strong moment condition. Take the sample mean of independent and identically distributed (i.i.d.) r.v.'s, for example. It is known that asymptotic normality holds under the second moment condition, and that an $r$-term Edgeworth expansion is valid under the $(r+2)$th moment condition plus some smoothness condition (e.g., a nonlattice or the Cramér condition). However, for the saddlepoint approximation one needs the much stronger condition that the exponential moment exists around the origin. This certainly limits the applicability of saddlepoint approximations in practice.

In this paper we shall focus on the saddlepoint approximation of the Student's $t$-statistic. It is common knowledge that the Student's $t$-statistic plays a pivotal role in statistics and is the most widely used statistic in the inference of a population mean. Therefore, accurate approximations to its d.f.'s become particularly important. Toward this end, Daniels and Young (1991) derived a saddlepoint approximation for the Student's $t$-statistic. However, their conditions are far too strong to be useful in practice. They require that the exponential moment of the square of the underlying r.v.'s exists near the origin. In other words, the underlying tail probability of the r.v.'s needs to go to zero as fast as the Normal distribution does. This is indeed a very severe restriction and makes the approximation hardly useful in practice. Even the exponential distribution can not satisfy this condition.

One of the purposes of this paper is to investigate how to weaken the strong moment condition given in Daniels and Young (1991) in the saddlepoint approximation of the Student's $t$-statistic. One of the key findings of the paper is that this very strong exponential moment condition can be totally eliminated. This result is highly significant in statistical inference for two reasons:

1. First of all, it makes the saddlepoint approximation more widely applicable. It is known [Giné, Götze, and Mason (1997)] that the Student's $t$-statistic is asymptotically $N(0,1)$ if and only if the r.v. is in the domain of attraction of the Normal law and that it has an $r$-term ($r \geq 1$) Edgeworth expansion under the $(r+2)$th moment condition plus some smoothness condition (e.g., the nonlattice or Cramér condition) [Hall (1987)]. Both asymptotic normality and Edgeworth expansion will not hold under heavy tail distributions, such as the Cauchy distribution. By contrast, this paper shows that the saddlepoint approximation does not need any moment condition at all and, at the same time, it provides an extremely accurate approximation to the tail probability of the Student's $t$-statistic.

2. Second, the fact that no moment condition is required for the saddlepoint approximation shows that the Student's $t$-statistic can guard against possible heavy tail distributions. This confirms the folklore that the Student's $t$-statistic is very robust against possible outliers.

For these reasons, the *saddlepoint approximation should always be used in practice whenever possible.*

The layout of the paper is as follows. Section 2 presents the formulation of the problem. Some notation and a brief review are given in Section 3. The main result will be presented in Section 4. Some numerical studies are given in Section 5. The proofs are given in Section 6. All technical details are left to the Appendix.

**2. Formulation of the problem.** Let $\{X, X_n, n \geq 1\}$ be a sequence of i.i.d. nondegenerate r.v.'s with d.f. $F(x)$. Write

$$\overline{X} = \frac{1}{n}\sum_{j=1}^{n} X_j.$$

Now consider the Student's $t$-statistic

$$T_n := \sqrt{n}\,\overline{X}/S, \qquad \text{where } S^2 = (n-1)^{-1}\sum_{j=1}^{n}(X_j - \overline{X})^2 \text{ for } n \geq 2.$$

It is known that asymptotic normality of $T_n$ holds if and only if $X$ is in the domain of attraction of the Normal law [Giné, Götze and Mason (1997)], which implies that $E|X|^{2-\varepsilon} < \infty$ for any $\varepsilon > 0$. Hall (1987) showed that $T_n$ has an $r$-term $(r \geq 1)$ Edgeworth expansion under the $(r+2)$th moment condition plus some smoothness condition (e.g., nonlattice or Cramér condition). On the other hand, Daniels and Young (1991) derived a saddlepoint approximation of Lugannani and Rice's (1980) type for the tail probability of $T_n$ under the assumption that the joint moment generating function of $X$ and $X^2$ exists, that is,

$$M(s,t) = \exp\{K(s,t)\} = Ee^{sX+tX^2} < \infty \tag{2.1}$$

for $(s,t)^T$ in a neighborhood of the origin. However, condition (2.1) requires that the tail probability of the underlying d.f. drop to zero at least as fast as a Normal r.v. does. This is, indeed, a very restrictive requirement; for example, it is violated even for the Exponential distribution. This severely limits its applicability in statistical inference. The natural question is: "*Is it possible to weaken the strong exponential moment condition and, if so, how far can we go?*"

Note that $T_n$ is closely related to the so-called *self-normalized sum* defined by

$$\frac{S_n}{V_n} = \sqrt{n}\frac{\overline{X}}{\overline{V}_n},$$

where

$$S_n = \sum_{i=1}^{n} X_i,$$

$$V_n^2 = \sum_{i=1}^{n} X_i^2,$$

$$\overline{V}_n = \left(n^{-1}\sum_{i=1}^{n} X_i^2\right)^{1/2}.$$

To see this, we note the following identity:

$$T_n = \frac{S_n}{V_n}\left(\frac{n-1}{n-(S_n/V_n)^2}\right)^{1/2}.$$

It suffices to investigate the self-normalized sum, $S_n/V_n$, because of the following identity:

$$\{T_n \geq t\} = \left\{\frac{S_n}{V_n} \geq t\left(\frac{n}{n+t^2-1}\right)^{1/2}\right\}. \tag{2.2}$$

There has been a growing literature on the study of self-normalized sums in recent years. For instance, one can refer to Logan, Mallows, Rice and Shepp (1973) for weak convergence, to Griffin and Kuelbs (1989, 1991) for a self-normalized law of the iterated logarithm, to Giné, Götze and Mason (1997) for the necessary and sufficient condition for asymptotic normality, and to Wang and Jing (1999) for the exponential nonuniform Berry–Esseen bound under finite moment conditions. However, the work most relevant to the present paper is that by Shao (1997), who studied self-normalized large deviations. Among other results, Shao [(1997), Corollary 1.1] showed the following result.

THEOREM 2.1 [Shao (1997)]. *Assume that either $EX = 0$ or $EX^2 = \infty$. Then for $x > 0$,*

$$\lim_{n\to\infty} P\left(\frac{\overline{X}}{\overline{V}_n} \geq x\right)^{1/n} = \sup_{c\geq 0}\inf_{t\geq 0} E\exp(t(cX - x(X^2+c^2)/2)).$$

Since for any r.v. $X$, either $EX^2 < \infty$ or $EX^2 = \infty$, the assumption that $EX = 0$ is reasonable if $EX^2 < \infty$. In other words, the large deviation for the self-normalized sum in Theorem 2.1 holds *without* assuming any moment conditions [see Remark 1.1 of Shao (1997)]. By contrast, a strong condition (2.1) is needed to derive the saddlepoint approximation for the self-normalized sum by Daniels and Young's approach, as noted earlier. This begs the question whether one can completely eliminate the condition (2.1) in the saddlepoint approximation of the self-normalized sum. The answer to this question is in the affirmative, as is shown later in the paper.

**3. Notation and brief review.** In this section we shall introduce some notation that is used in later sections. We do this by briefly deriving saddlepoint approximations of the self-normalized sum $S_n/V_n$ under the strong exponential moment condition (2.1), following similar lines to those in Daniels and Young (1991).

The first step involves finding the saddlepoint approximations of the joint density of $(\overline{X}, \overline{Y})^T$, where $Y = X^2$, $Y_i = X_i^2$ for $1 \leq i \leq n$ and $\overline{Y} =$

$n^{-1}\sum_{i=1}^n Y_i$. Assume that the cumulant-generating function of $(X, X^2)^T$ satisfies

$$K(s,t) = \ln M(s,t) = \ln Ee^{sX+tX^2} < \infty \tag{3.1}$$

in a neighborhood of the origin. Denote

$$\begin{aligned} K_s(s,t) &:= \frac{\partial K(s,t)}{\partial s}, \\ K_t(s,t) &:= \frac{\partial K(s,t)}{\partial t}, \\ K_{ss}(s,t) &:= \frac{\partial^2 K(s,t)}{\partial s^2} \qquad \text{and so on.} \end{aligned}$$

Assume that $(X, X^2)^T$ has an integrable characteristic function. Then, by the Fourier inversion formula, the saddlepoint approximation to the joint density, $f_n(x,y)$, of $(\overline{X}, \overline{Y})^T$ is given by

$$f_n(x,y) = \frac{n^2}{(2\pi i)^2} \iint e^{-n[sx+ty-K(s,t)]}\, ds\, dt = \hat{f}_n(x,y)(1 + r_n/n), \tag{3.2}$$

where integration is along admissible paths in $R^2$, and

$$\hat{f}_n(x,y) = \frac{n}{2\pi} \frac{e^{-n[\hat{s}x+\hat{t}y-K(\hat{s},\hat{t})]}}{[K_{ss}(\hat{s},\hat{t})K_{tt}(\hat{s},\hat{t}) - K_{st}^2(\hat{s},\hat{t})]^{1/2}}, \tag{3.3}$$

where $\hat{s} = \hat{s}(x,y)$ and $\hat{t} = \hat{t}(x,y)$ are solutions to

$$K_s(\hat{s},\hat{t}) = x, \qquad K_t(\hat{s},\hat{t}) = y, \tag{3.4}$$

and $|r_n| < C$ for some $C > 0$ if $(x,y)^T$ is contained in a compact set.

The second step is to find the joint density $f_{(\overline{X}, \overline{X}/\overline{V}_n)}(a,b)$. Let $a = x$, $b = x/\sqrt{y}$ $(y > 0)$. The inverse transformation and its Jacobian determinant are

$$x \equiv x(a,b) := a, \qquad y \equiv y(a,b) := a^2/b^2, \qquad J(a,b) = 2a^2/b^3. \tag{3.5}$$

Thus, the saddlepoint approximation to the joint density of $(\overline{X}, \overline{X}/\overline{V}_n)^T$ is

$$\hat{f}_{(\overline{X}, \overline{X}/\overline{V}_n)}(a,b) = |J(a,b)|\hat{f}_n(x,y) = \frac{n}{2\pi} \frac{|J(a,b)|}{\det\{\Delta(a,b)\}^{1/2}} e^{-n\Lambda(a,b)},$$

where $\hat{s} = \hat{s}(x(a,b), y(a,b))$, $\hat{t} = \hat{t}(x(a,b), y(a,b))$, and

$$\begin{aligned} \Lambda(a,b) &= \hat{s}x(a,b) + \hat{t}y(a,b) - K(\hat{s},\hat{t}) = \hat{s}a + \hat{t}a^2/b^2 - K(\hat{s},\hat{t}), \\ \Delta(a,b) &= \begin{pmatrix} K_{ss}(\hat{s},\hat{t}) & K_{st}(\hat{s},\hat{t}) \\ K_{st}(\hat{s},\hat{t}) & K_{tt}(\hat{s},\hat{t}) \end{pmatrix}, \end{aligned}$$

where $\hat{s}$ and $\hat{t}$ satisfy $K_s(\hat{s},\hat{t}) = a$ and $K_t(\hat{s},\hat{t}) = a^2/b^2$. After some simple algebra, we obtain

$$\Lambda_a(a,b) = \hat{s} + \frac{2\hat{t}a}{b^2},$$

$$\Lambda_b(a,b) = -\frac{2\hat{t}a^2}{b^3},$$

$$\Lambda_{aa}(a,b) = \frac{2\hat{t}}{b^2} + \left(1, \frac{2a}{b^2}\right)\Delta(a,b)^{-1}\left(1, \frac{2a}{b^2}\right)^T.$$

The third step involves finding the marginal density of $\overline{X}/\overline{V}_n$. Let $a_0 = a_0(b)$ be such that

$$\Lambda(a_0,b) := \inf_a \Lambda_a(a,b).$$

If we assume that $\Lambda_{aa}(a_0,b) > 0$, then $a_0 = a_0(b)$ is the unique solution of $\Lambda_a(a_0,b) = 0$. Then the Laplace approximation of the marginal density of $\overline{X}/\overline{V}_n$ is

$$\hat{f}_{\overline{X}/\overline{V}_n}(b) =: \sqrt{\frac{n}{2\pi}}\frac{|J(a_0,b)|}{\det\{\Delta(a_0,b)\}^{1/2}\Lambda_{aa}^{1/2}(a_0,b)}e^{-n\Lambda(a_0,b)}.$$

Finally, by applying another Laplace approximation in integrating $\hat{f}_{(\overline{X}/\overline{V}_n)}(b)$, we get the saddlepoint approximation for the self-normalized sum. We summarize the result in the following theorem.

THEOREM 3.1. *Assume that:*

(C1) $Ee^{i\xi X + i\eta X^2} \in L^v(R^2)$ *for some* $v > 1$, *that is,* $\int\int |Ee^{i\xi X + i\eta X^2}|^v \, d\xi \, d\eta < \infty$.

(C2) $\Lambda_{aa}(a_0,b) > 0$.

(C3) $Ee^{sX+tX^2} < \infty$ *in a neighborhood of the origin.*

*Then we have*

$$P\left(\frac{\overline{X}}{\overline{V}_n} \geq b\right) = 1 - \Phi(\sqrt{n}w) - \frac{\phi(\sqrt{n}w)}{\sqrt{n}}\left(\frac{1}{w} - \frac{1}{v} + O(n^{-1})\right),$$

*where* $w = \sqrt{2\Lambda(a_0,b)}$ *and* $v = -\det\{\Delta(a_0,b)\}^{1/2}\Lambda_{aa}^{1/2}(a_0,b)\hat{t}_0$, *and* $(\hat{s}_0, \hat{t}_0, a_0)$ *are solutions* $(s,t,a)$ *to the equations*

$$s + \frac{2ta}{b^2} = 0, \qquad \frac{EXe^{sX+tX^2}}{Ee^{sX+tX^2}} = a, \qquad \frac{EX^2e^{sX+tX^2}}{Ee^{sX+tX^2}} = \frac{a^2}{b^2}. \tag{3.6}$$

REMARK 3.1. From the first equation of (3.6), we obtain that $s = -2ta/b^2$. Therefore, on substituting this into the other two equations, then (3.6) reduces to

$$\frac{EXe^{t(-2aX/b^2+X^2)}}{Ee^{t(-2aX/b^2+X^2)}} = a, \qquad \frac{EX^2e^{t(-2aX/b^2+X^2)}}{Ee^{t(-2aX/b^2+X^2)}} = \frac{a^2}{b^2}. \tag{3.7}$$

**4. Main results.** The saddlepoint approximation for the Student's $t$-statistic under the strong exponential moment conditions was given in Theorem 3.1. Condition (C1) is a smoothness condition, which validates the Fourier inversion formula (3.2). It is satisfied, for instance, when the r.v. $X$ has a density function. The main purpose of this paper is to remove conditions (C2) and (C3) in Theorem 3.1.

THEOREM 4.1. *Let $0 < b < 1$ and let $X$ be a r.v. with $EX = 0$ or $EX^2 = \infty$. Assume further that condition* (C1) *in Theorem* 3.1 *holds. Then*

$$P\left(\frac{\overline{X}}{\overline{V}_n} \geq b\right) = 1 - \Phi(\sqrt{n}w) - \frac{\phi(\sqrt{n}w)}{\sqrt{n}}\left(\frac{1}{w} - \frac{1}{v} + O(n^{-1})\right), \tag{4.1}$$

*where $w$ and $v$ are defined the same as in Theorem* 3.1.

We make the following remarks:

1. When $-1 < b < 0$, similarly, we have
$$P(\overline{X}/\overline{V}_n \leq b) = \Phi(\sqrt{n}w) - \frac{\phi(\sqrt{n}w)}{\sqrt{n}}\left(-\frac{1}{w} - \frac{1}{v} + O(n^{-1})\right).$$
2. Theorem 4.1 remains valid when $b = \pm 1$. Take $b = 1$ for instance. From Proposition 6.1, condition (C1) implies that $X$ is a continuous r.v. Then the left-hand side of (4.1) is
$$P(\overline{X}/\overline{V}_n \geq b) = P(X_1 = \cdots = X_n, X_1 > 0) = 0.$$
On the other hand, it can be shown that $w = \infty$ if $b = 1$ (see Remark **??**), which implies that the right-hand side of (4.1) is also zero.
3. The case for $b = 0$ is slightly different. By the Berry–Esseen bound, we have
$$P(\overline{X}/\overline{V}_n \geq 0) = P(\overline{X} \geq 0) = \tfrac{1}{2}\{1 + O(n^{-1/2})\},$$
provided that $E|X|^3 < \infty$, which is the minimal moment condition required here. Comparing this with Theorem 4.1, we notice that a stronger condition is needed for the case when $b = 0$ than when $b \neq 0$. It may seem odd that one needs stronger conditions in the middle of the distribution than in the tails. The reason is that when $b = 0$ there is nothing to offset the effect of possibly heavy tail distributions. Therefore, one must impose extra conditions to control the tail behavior.

**5. Numerical study.** In this section we conduct some numerical studies to investigate the performance of the saddlepoint approximation for the Student's $t$-statistic. Let $X_1, \ldots, X_n$ be a random sample from a distribution with p.d.f. $f(x)$. We shall choose $f(x)$ from several well-known density functions, ranging from one with very thin tails (e.g., Normal density) to one with rather heavy tails (e.g., Cauchy).

Our interest is to calculate the probability of the self-normalized sum, $P(\overline{X}/\overline{V}_n \geq b)$, for a range of values of $b \in (0,1)$. Since the exact value of the above probability is difficult to obtain in practice, we calculate its "exact" probability by 1,000,000 Monte Carlo simulations. Then, we compare how well the saddlepoint approximation performs in comparison with other approximation methods, such as the large deviation [Shao (1997)], the Edgeworth expansion [Hall (1987)], and the Normal approximation.

For illustration purposes, we choose the sample size to be $n = 5$, since different sample sizes display similar patterns. In the tables below, we use the following abbreviations:

"True" = true probability,
"Saddle" = saddlepoint approximation,
"Edgeworth" = Edgeworth expansion,
"L.D." = large deviation,
"N.A." = Normal approximation,
"R.E." = relative error.

5.1. *Saddlepoint approximation vs. large deviation.* Here we compare the saddlepoint approximation for self-normalized sums and large deviation results of Shao (1997). Let $X_1, \ldots, X_n$ be a random sample from the Standard Normal distribution with p.d.f.

$$f(x) = \frac{1}{\sqrt{2\pi}} e^{-x^2/2}.$$

The reason for deliberately choosing this "nicest" density function is based on the belief that any approximation method should probably work at its "best" under this special situation if it works at all. In other words, if a method does not work well in this case, we cannot expect it to work well in other cases either. The simulation results are presented in Table 1 and Figure 1.

We first make some general remarks.

(i) First of all, the saddlepoint approximation provides extremely accurate approximations to the exact probabilities and performs uniformly better than the other approximation methods, even for sample sizes as small as 5. In fact, Figure 1 shows that the saddlepoint approximation is almost indistinguishable from the true probability. The superiority of the saddlepoint approximation becomes even more pronounced in the tails of the distributions.

(ii) Since the sample is from a Normal distribution, the Normal approximation and one-term Edgeworth expansion to $P(\overline{X}/\overline{V}_n \geq b)$ coincide. Table 1 shows that the Normal approximation gives very good approximation at the center of the distribution in this "nicest" case. However, the approximation soon starts to deteriorate very quickly toward the tail area of the distribution.

(iii) The large deviation performs miserably throughout the whole range. It is much worse than even the Normal approximation at the center of the distribution. In the tail area, the saddlepoint approximation is much superior to the large deviation. This shows that one can NOT rely on the large deviation to give accurate approximations of probabilities.

This example clearly demonstrates that the large deviation is no substitute for the saddlepoint approximation when it comes to accurate approximations, even for a case as nice as the Normal distribution. The same phenomenon has also been found for other underlying d.f.'s. For this reason, we shall not include the large deviation in our simulation studies below.

To see why the large deviation performs so poorly, we note that Theorem 2.1 gives the limit of $P(\overline{X}/\overline{V}_n \geq b)^{1/n}$ as $n \to \infty$. However, $[C_n P(\overline{X}/\overline{V}_n \geq$ $b)]^{1/n}$ would give the same limit as long as $C_n^{1/n} \to 1$. That is, the large deviation only captures the exponential component and any other terms are simply thrown away.

TABLE 1
$f(x) = (2\pi)^{-1/2} e^{-x^2/2}$ *(Normal density)*

| *b* | **True** | **Saddle** | **(R.E.)** | **L.D.** | **(R.E.)** | **N.A.** | **(R.E.)** |
|---|---|---|---|---|---|---|---|
| 0.05 | 0.4621 | 0.4621 | (0.0001) | 0.9938 | (1.15) | 0.4555 | (0.01) |
| 0.10 | 0.4243 | 0.4244 | (0.0003) | 0.9752 | (1.30) | 0.4115 | (0.03) |
| 0.15 | 0.3869 | 0.3872 | (0.0007) | 0.9447 | (1.44) | 0.3687 | (0.05) |
| 0.20 | 0.3500 | 0.3505 | (0.001) | 0.9030 | (1.58) | 0.3274 | (0.06) |
| 0.25 | 0.3138 | 0.3146 | (0.003) | 0.8510 | (1.71) | 0.2881 | (0.08) |
| 0.30 | 0.2785 | 0.2797 | (0.004) | 0.7900 | (1.84) | 0.2512 | (0.10) |
| 0.35 | 0.2443 | 0.2460 | (0.007) | 0.7213 | (1.95) | 0.2169 | (0.11) |
| 0.40 | 0.2113 | 0.2136 | (0.01) | 0.6467 | (2.06) | 0.1855 | (0.12) |
| 0.45 | 0.1799 | 0.1829 | (0.02) | 0.5680 | (2.16) | 0.1572 | (0.13) |
| 0.50 | 0.1502 | 0.1539 | (0.02) | 0.4871 | (2.24) | 0.1318 | (0.12) |
| 0.55 | 0.1225 | 0.1268 | (0.04) | 0.4063 | (2.32) | 0.1094 | (0.11) |
| 0.60 | 0.0970 | 0.1019 | (0.05) | 0.3277 | (2.38) | 0.0899 | (0.07) |
| 0.65 | 0.0739 | 0.0793 | (0.07) | 0.2534 | (2.42) | 0.0731 | (0.01) |
| 0.70 | 0.0536 | 0.0592 | (0.10) | 0.1857 | (2.46) | 0.0588 | (0.10) |
| 0.75 | 0.0363 | 0.0417 | (0.15) | 0.1266 | (2.49) | 0.0468 | (0.29) |
| 0.80 | 0.0223 | 0.0271 | (0.22) | 0.0778 | (2.49) | 0.0368 | (0.65) |
| 0.85 | 0.0116 | 0.0154 | (0.33) | 0.0406 | (2.49) | 0.0287 | (1.46) |
| 0.90 | 0.0045 | 0.0070 | (0.53) | 0.0157 | (2.46) | 0.0221 | (3.86) |
| 0.95 | 0.0009 | 0.0018 | (1.03) | 0.0030 | (2.42) | 0.0168 | (18.4) |

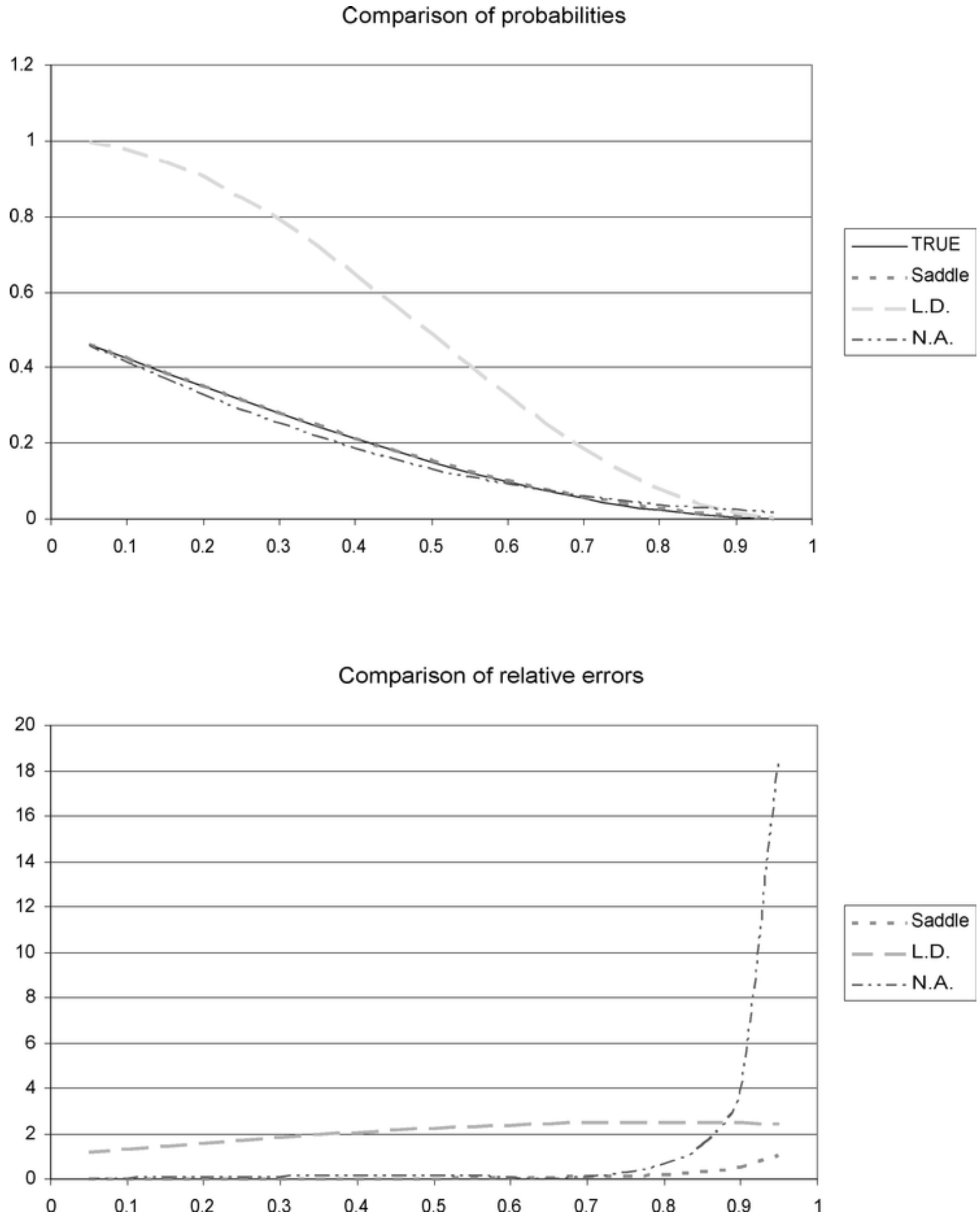

FIG. 1. *Comparisons under the Normal density.*

In a way, the relationship between the large deviation and the saddlepoint approximation is a little like that between the Normal approximation and the Edgeworth expansion, since in both cases, the former provides the dominant term for the latter. One major difference is the following. The Normal approximation can be used in statistical inference when the sample size is reasonably large and the Edgeworth expansion can often provide more accurate approximations than the Normal approximation. However, one can not usually rely on large deviation probability to calculate tail probabilities in general since the approximations are often too crude to be useful, as shown in the last example. By contrast, the saddlepoint method can provide extremely accurate approximations throughout the range.

5.2. *Saddlepoint approximations for light tailed distributions.* Here, we study the accuracy of the saddlepoint approximation to $P(\overline{X}/\overline{V}_n \geq x)$ when the underlying distribution has thin tails. Let $X_1, \ldots, X_n$ be a random sample from the centered exponential density with p.d.f.

$$f(x) = e^{-(x+1)}, \qquad x \geq -1.$$

The tail of the density decreases exponentially fast (but not as fast as the Normal density function). As mentioned before, even for this "nice" density, the stringent exponential moment condition given by Daniels and Young (1991) is not satisfied. But the saddlepoint approximation still holds from Theorem 4.1. The Normal approximation and the Edgeworth expansion are included for comparison. The results are presented in Table 2 and Figure 2.

We make the following observations.

(i) The saddlepoint approximation is remarkably accurate and uniformly better than the other approximation methods. Most of the relative errors fall below 10%, and the maximum error is only 17% near the center of the distribution.

(ii) The Edgeworth expansion performs better than the Normal approximation throughout the whole range. Both give reasonable approximations at the center, but they turn very bad toward the tail areas, where the relative errors are of the order of 1000% for tail area probabilities in the order of 1%. By comparison, the errors for the saddlepoint approximation do not exceed 20% for the whole region.

(iii) This example clearly demonstrates why accurate approximations of the tail area probabilities are important in statistical inference. It is easy to conceive of a hypothesis test such that its $p$-value is given by $P_{H_0}(\overline{X}/\overline{V}_n \geq 0.75)$, where the $X_i$'s follow a centered exponential distribution under $H_0$. From Table 2, the true value is $0.0088 < 0.01$, which leads to the rejection of $H_0$ at significance level 1%. The same conclusion would be reached by using the saddlepoint approximation, but not by using the Normal approximation or the Edgeworth expansion.

5.3. *Saddlepoint approximations for heavy tailed distributions.* Here we are interested in the accuracy of the saddlepoint approximation for self-normalized sums when the underlying distribution has heavy tails. We shall give two examples.

EXAMPLE 5.1. Let $X_1, \ldots, X_n$ be a random sample from the $t_2$ distribution with p.d.f.

$$f(x) = \frac{1}{2^{3/2}(1 + x^2/2)^{3/2}}.$$

Clearly, $EX_1 = 0$ and $\mathrm{Var}(X_1) = \infty$. Also, it is easy to check that $X_1$ is in the domain of attraction of the Normal law. It then follows from Giné, Götze and Mason (1997) that the Student's $t$-statistic is asymptotically $N(0,1)$. Clearly, the saddlepoint approximation still holds under this heavy tail distribution, following Theorem 4.1. So, in this case, we can compare the saddlepoint approximation with the Normal approximation of the Student's $t$-statistic. The results are summarized in Table 3 and Figure 3.

TABLE 2
$f(x) = e^{-(x+1)}$, $x \geq -1$ *(centered exponential density)*

| $b$ | **True** | **Saddle** | **(R.E.)** | **Normal** | **(R.E.)** | **Edgeworth** | **(R.E.)** |
|---|---|---|---|---|---|---|---|
| 0.05 | 0.4231 | 0.4951 | (0.170) | 0.4602 | (0.09) | 0.4024 | (0.05) |
| 0.10 | 0.3869 | 0.4267 | (0.103) | 0.4207 | (0.09) | 0.3611 | (0.07) |
| 0.15 | 0.3487 | 0.3486 | (0.000) | 0.3821 | (0.10) | 0.3197 | (0.08) |
| 0.20 | 0.3090 | 0.3046 | (0.001) | 0.3446 | (0.12) | 0.2792 | (0.10) |
| 0.25 | 0.2680 | 0.2633 | (0.018) | 0.3085 | (0.15) | 0.2407 | (0.10) |
| 0.30 | 0.2270 | 0.2223 | (0.021) | 0.2743 | (0.21) | 0.2052 | (0.10) |
| 0.35 | 0.1866 | 0.1825 | (0.022) | 0.2420 | (0.20) | 0.1732 | (0.07) |
| 0.40 | 0.1486 | 0.1451 | (0.023) | 0.2119 | (0.43) | 0.1452 | (0.02) |
| 0.45 | 0.1141 | 0.1114 | (0.024) | 0.1841 | (0.61) | 0.1214 | (0.06) |
| 0.50 | 0.0840 | 0.0822 | (0.022) | 0.1587 | (0.89) | 0.1015 | (0.21) |
| 0.55 | 0.0594 | 0.0581 | (0.023) | 0.1357 | (1.28) | 0.0851 | (0.43) |
| 0.60 | 0.0402 | 0.0391 | (0.028) | 0.1151 | (1.86) | 0.0717 | (0.78) |
| 0.65 | 0.0256 | 0.0250 | (0.026) | 0.0968 | (2.77) | 0.0608 | (1.37) |
| 0.70 | 0.0156 | 0.0151 | (0.033) | 0.0808 | (4.19) | 0.0517 | (2.32) |
| 0.75 | 0.0088 | 0.0085 | (0.039) | 0.0668 | (6.59) | 0.0441 | (4.01) |
| 0.80 | 0.0045 | 0.0044 | (0.029) | 0.0548 | (11.15) | 0.0376 | (7.34) |
| 0.85 | 0.0021 | 0.0020 | (0.031) | 0.0446 | (20.43) | 0.0319 | (14.34) |
| 0.90 | 0.0008 | 0.00075 | (0.064) | 0.0359 | (43.86) | 0.0269 | (32.57) |

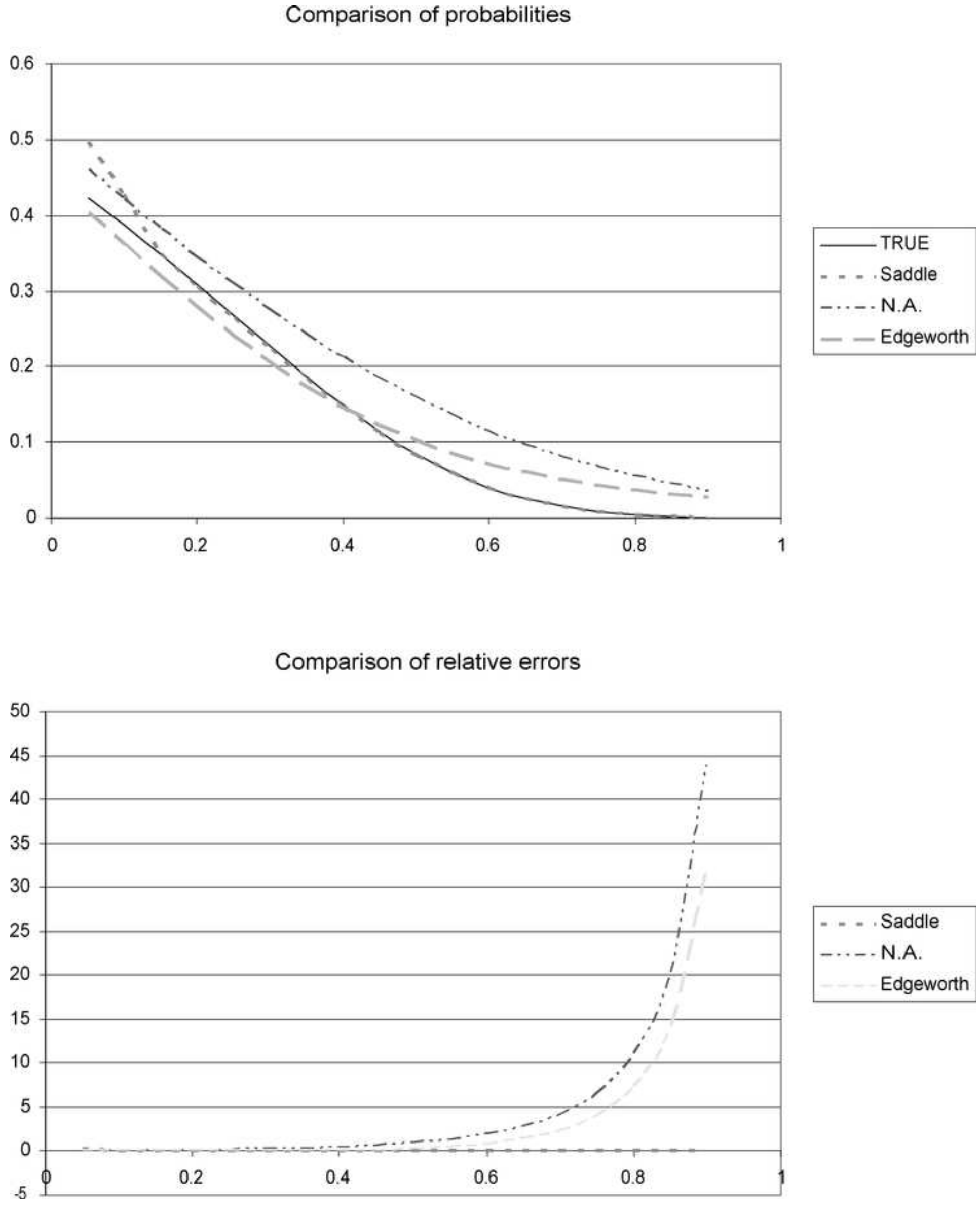

Fig. 2. *Comparisons under exponential density.*

Example 5.2. Let $X_1, \dots, X_n$ be a random sample from the Cauchy distribution with p.d.f.

$$f(x) = \frac{1}{\pi(1+x^2)}.$$

Note that the usual Normal approximation and Edgeworth expansion do not exist here. However, the saddlepoint approximation continues to hold here. The results are given in Table 4 and Figure 4.

We make some remarks about the two examples.

(i) Clearly, the saddlepoint approximation is remarkably accurate even for these rather heavy tail distributions. The relative errors remain very small (under 11% and 13%, resp.) for the range considered.

(ii) For the $t_2$ density case, asymptotic normality holds and the Normal approximation performs rather well in the center, but it becomes very poor toward the tail area. In fact, the relative errors start to shoot up just as the tail probability decreases to around 5% and beyond, which is the area of most interest in statistical inference. The plot of relative errors in Figure 3 should leave all doubts behind.

TABLE 3
$f(x) = 2^{-3/2}(1 + x^2/2)^{-3/2}$ *($t_2$ density)*

| $b$ | **True** | **Saddle** | **(R.E.)** | **N.A.** | **(R.E.)** |
|---|---|---|---|---|---|
| 0.40 | 0.2386 | 0.2637 | (0.105) | 0.2119 | (0.11) |
| 0.45 | 0.1987 | 0.2146 | (0.080) | 0.1841 | (0.07) |
| 0.50 | 0.1598 | 0.1708 | (0.069) | 0.1587 | (0.01) |
| 0.55 | 0.1255 | 0.1322 | (0.053) | 0.1357 | (0.08) |
| 0.60 | 0.0953 | 0.0990 | (0.040) | 0.1151 | (0.21) |
| 0.65 | 0.0694 | 0.0713 | (0.027) | 0.0968 | (0.39) |
| 0.70 | 0.0479 | 0.0488 | (0.019) | 0.0808 | (0.69) |
| 0.75 | 0.0310 | 0.0312 | (0.007) | 0.0668 | (1.15) |
| 0.80 | 0.0183 | 0.0182 | (0.006) | 0.0548 | (2.00) |
| 0.85 | 0.0094 | 0.0093 | (0.019) | 0.0446 | (3.72) |
| 0.90 | 0.0038 | 0.0036 | (0.056) | 0.0359 | (8.34) |

TABLE 4
$f(x) = \pi^{-1}(1 + x^2)^{-1}$ *(Cauchy density)*

| $b$ | **True** | **Saddle** | **(R.E.)** |
|---|---|---|---|
| 0.40 | 0.2712 | 0.3058 | (0.13) |
| 0.45 | 0.2085 | 0.2302 | (0.10) |
| 0.50 | 0.1515 | 0.1697 | (0.12) |
| 0.55 | 0.1117 | 0.1218 | (0.09) |
| 0.60 | 0.0798 | 0.0845 | (0.06) |
| 0.65 | 0.0537 | 0.0563 | (0.05) |
| 0.70 | 0.0344 | 0.0356 | (0.04) |
| 0.75 | 0.0207 | 0.0210 | (0.02) |
| 0.80 | 0.0112 | 0.0113 | (0.01) |
| 0.85 | 0.0052 | 0.0052 | (0.00) |
| 0.90 | 0.0019 | 0.0019 | (0.02) |

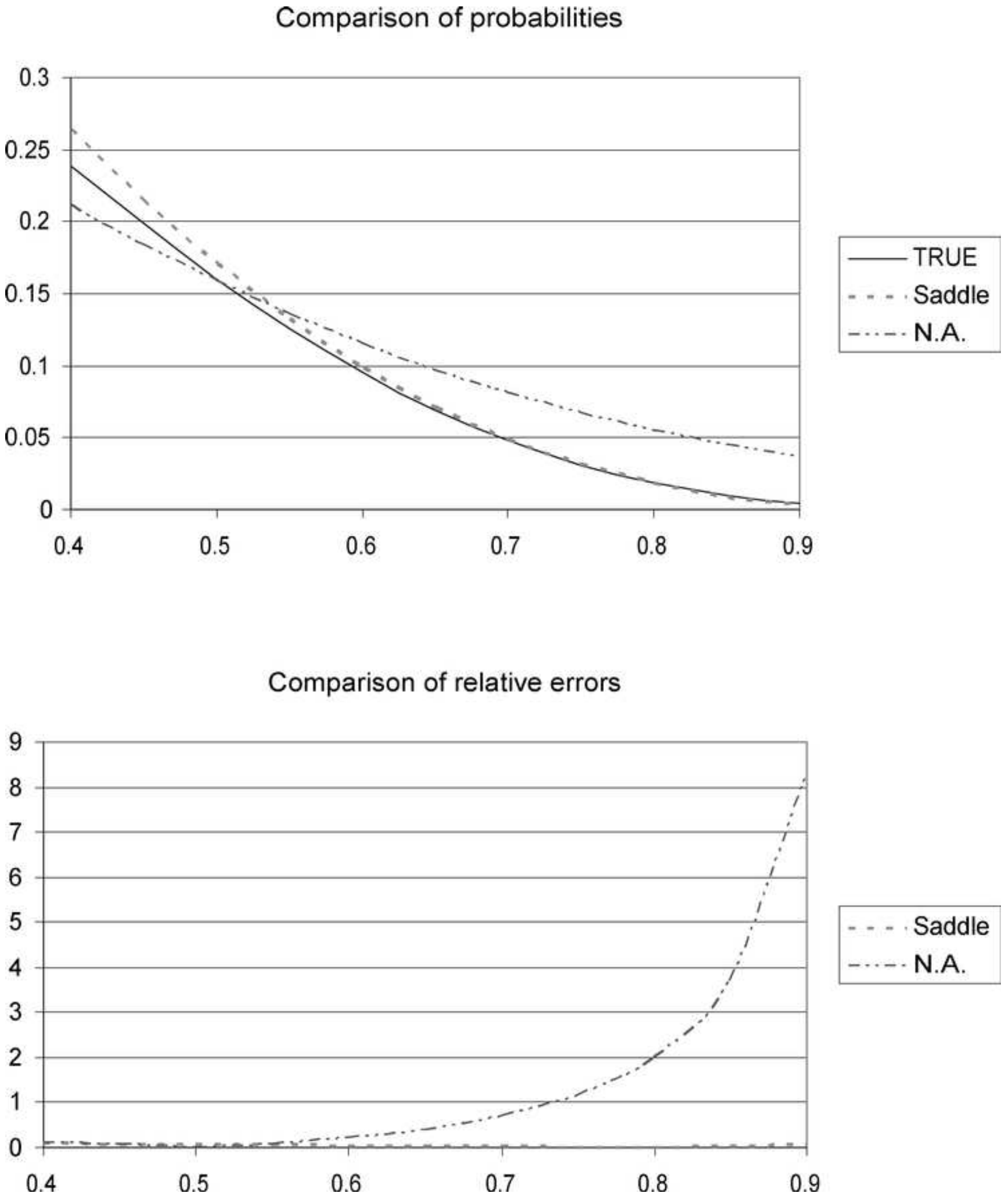

FIG. 3. *Comparisons under the $t_2$ density.*

(iii) We have seen that the saddlepoint approximation provides extremely accurate approximation of the distribution of the self-normalized sum or, equivalently, of the Student's $t$-statistic, particularly near the tail area. It is also clear that the tail probability of the Student's $t$-statistic decreases exponentially fast. These properties hold irrespective of whether the underlying density has light or heavy tails. These results confirm the common belief that the Student's $t$-statistic provides a very *robust* procedure for the statistical inference of a population mean with a possible heavy-tailed dis-tribution. On the other hand, it is well known that the sample mean is very sensitive to outliers and is not robust against heavy-tailed distributions.

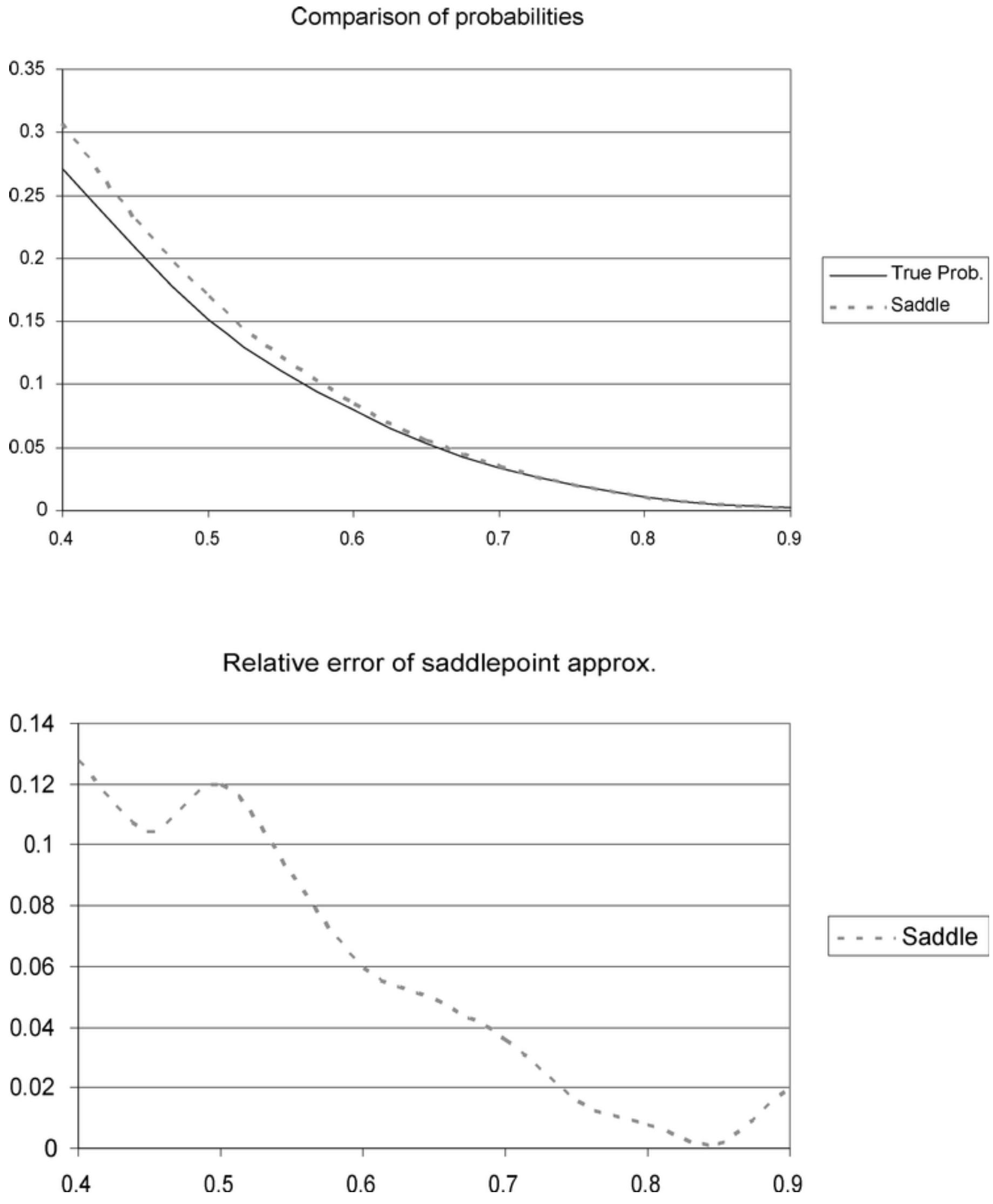

FIG. 4. *Comparisons under the Cauchy density.*

(iv) Robustness of the self-normalized sums or, equivalently, the Student's $t$-statistic, can also be explained intuitively as follows. It is well known that when there is an outlier on the right-hand side among the observations $X_1, \ldots, X_n$, the sample mean, $\overline{X}$, is dominated by the largest order statis-tic, $X_{(n)} = \max\{X_1, \dots, X_n\}$. For self-normalized sums, $\overline{X}/\overline{V}_n$, both $\overline{X}$ and $\overline{V}_n$ are dominated by $X_{(n)}$, effectively cancelling the influence of any outlier.

5.4. *Summary.* The Student's $t$-statistic is one of the most commonly used statistics in inference. We have derived a saddlepoint approximation for the Student's $t$-statistic under no moment condition. The key results are summarized as follows.

1. The saddlepoint approximation provides extremely accurate approximations to the distribution of the Student's $t$-statistic. The approximation is particularly useful in calculating small probabilities in the tail areas, which are often of great interest in practice.
2. The saddlepoint approximation holds under no moment condition. This makes the application of the saddlepoint approximation very broad. This is significant for the user since one can use the approximation without having to worry about whether or not the result is valid.
3. The Student's $t$-statistic is very robust against possible outliers.

For those reasons, the saddlepoint approximation of the Student's $t$-statistic should always be used in practice whenever possible.

## 6. Proof of Theorem 4.1.

The immediate consequence of condition (C1) is as follows.

PROPOSITION 6.1. *$F(x)$ is a continuous d.f. under condition* (C1) *of Theorem* 3.1.

PROOF. Let $2u$ be the smallest even integer not less than $v$. Then

$$\begin{aligned} &Ee^{i\xi(X_1+\cdots+X_u-X_{u+1}-\cdots-X_{2u})+i\eta(X_1^2+\cdots+X_u^2-X_{u+1}^2-\cdots-X_{2u}^2)} \\ &\quad = |Ee^{i\xi X_1+i\eta X_1^2}|^{2u} \in L^1(R^2). \end{aligned}$$

By the Fourier inversion theorem in $R^2$ [e.g., see (7.14) of Feller (1971)], $(X_1+\cdots+X_u-X_{u+1}-\cdots-X_{2u}, X_1^2+\cdots+X_u^2-X_{u+1}^2-\cdots-X_{2u}^2)^T$ has a bounded continuous density, which implies that $F(x)$ is a continuous d.f. $\square$

The key to getting rid of condition (C2) is the following.

PROPOSITION 6.2. *Assume that $F(x)$ is a continuous d.f. Then for each fixed $b \in (0,1)$, $\inf_{a>0} \Lambda(a,b)$ is attained at some finite unique point, $a_0 := a_0(b)$, which is the solution to $\Lambda_a(a,b) = 0$.*

PROOF. The proof follows from Lemmas A.6 and A.7. □

In an effort to remove condition (C3), we shall give the following two propositions.

PROPOSITION 6.3. *Under the conditions of Theorem 2.1, we have*

$$\begin{aligned}\lim_{n\to\infty} P(\overline{X}/\overline{V}_n \geq b)^{1/n} &= \sup_{a\geq 0}\inf_{t\geq 0} E\exp(t(aX - b(X^2+a^2)/2))\\ &= \exp\biggl\{-\inf_{a>0}\Lambda(a,b)\biggr\}.\end{aligned}$$

PROOF. The first equality follows from Theorem 1.1 of Shao (1997). The second one follows since

$$\begin{aligned}&\log\Bigl(\sup_{a\geq 0}\inf_{t\geq 0} E\exp\{t(aX - b(X^2+a^2)/2)\}\Bigr)\\ &\quad= -\inf_{a\geq 0}\sup_{t\geq 0}\Bigl(\frac{1}{2}tba^2 - \log E\exp\{t(aX - bX^2/2)\}\Bigr)\\ &\quad= -\inf_{a\geq 0}\sup_{t_1\leq 0}\Bigl(-t_1a^2 - \log E\exp\Bigl\{t_1\Bigl(-\frac{2a}{b}X + X^2\Bigr)\Bigr\}\Bigr)\\ &\qquad\qquad\qquad\qquad\qquad\qquad\qquad\qquad\qquad (\text{where } t_1 = -tb/2)\\ &\quad= -\inf_{a_1\geq 0}\sup_{t_1\leq 0}\Bigl(-t_1\frac{a_1^2}{b^2} - K\Bigl(-\frac{2a_1}{b^2}t_1, t_1\Bigr)\Bigr) \qquad (\text{where } a = a_1/b)\\ &\quad= -\inf_{a>0}\Lambda(a,b), \qquad\qquad [\text{by (A.2) and Lemma A.7}].\end{aligned}$$

□

The proposition establishes the relationship between the saddlepoint approximation formula of Theorem 3.1 and the large deviation results of Theorem 2.1. It shows that the dominant term in the saddlepoint approximation given in Theorem 3.1 is the same as that in the large deviation of Shao (1997). Since the latter requires no moment conditions at all, it is therefore reasonable to expect that Theorem 3.1 holds under no moment conditions as well. Unfortunately, the techniques used in Shao (1997) cannot be employed here for our purposes. One crucial result is the following.

PROPOSITION 6.4. *Assume that $F(x)$ is a continuous d.f. Then, for $0 < b < 1$, there exist solutions $(\hat{s}_0, \hat{t}_0, a_0)$ in (3.6) such that $\hat{s}_0 > 0$, $\hat{t}_0 < 0$ and $a_0 > 0$.*

PROOF. The proof follows straightaway from Lemmas A.3, A.6 and Remark 3.1. $\square$

The critical observation here is that $\hat{t}_0 < 0$, which implies that the cumulant generating function, $K(s,t) = \ln Ee^{sX+tX^2}$, always exists for $(s,t)^T$ in a small neighborhood of $(\hat{s}_0, \hat{t}_0)^T$ by the continuity of $K(s,t)$. This suggests that, in order to derive self-normalized saddlepoint approximations without moment conditions, we need to divide the probability, $P(\overline{X}/\overline{V}_n \geq b)$, into two regions:

(i) a small neighborhood of $(\hat{s}_0, \hat{t}_0)^T$ for which $\hat{t}_0 < 0$, where we need to show that there exists a saddlepoint approximation without any moment conditions;

(ii) the remaining region outside this small neighborhood of $(\hat{s}_0, \hat{t}_0)^T$, where we need to show that the probability is "negligible."

To make these statements precise, define

$$\begin{aligned}
\Omega(b) &= \{(x,y)^T | b \leq x/\sqrt{y} \leq 1\},\\
\Omega_0(b) &= \{(x,y)^T | (x-a_0)^2 + (y - a_0^2/b^2)^2 \leq \varepsilon^2\} \cap \Omega(b),\\
\Omega_1(b) &= \Omega(b) \setminus \Omega_0(b).
\end{aligned}$$

The closure of an arbitrary set, $A$, will be denoted as $A^-$. The plots of these regions are illustrated in Figure 5.

Hence, for any $0 < b < 1$,

$$\begin{aligned}
P(\overline{X} \geq b\overline{V}_n) &= \iint_{\Omega(b)} f_{(\overline{X},\overline{Y})}(x,y)\,dx\,dy\\
&= \iint_{\Omega_0(b)} f_{(\overline{X},\overline{Y})}(x,y)\,dx\,dy + P((\overline{X},\overline{Y})^T \in \Omega_1(b)) \qquad (6.1)\\
&:= J_1(b) + J_2(b).
\end{aligned}$$

Thus, the proof of Theorem 4.1 follows from the next two propositions.

PROPOSITION 6.5. *Under the conditions of Theorem 4.1 we have*

$$J_1(b) = 1 - \Phi(\sqrt{n}w) - \frac{\phi(\sqrt{n}w)}{\sqrt{n}}\left(\frac{1}{w} - \frac{1}{v} + O(n^{-1})\right), \qquad (6.2)$$

*where $w$ and $v$ are defined the same as in Theorem 3.1.*

PROOF. Denote $h_1(s,t;a,b) = K_s(s,t) - a$, $h_2(s,t;a,b) = K_t(s,t) - a^2/b^2$. Since $h_1(-\frac{2a_0}{b^2}\hat{t}_0, \hat{t}_0; a_0, b) = 0$, $h_2(-\frac{2a_0}{b^2}\hat{t}_0, \hat{t}_0; a_0, b) = 0$ and

$$\left.\begin{pmatrix} \dfrac{\partial h_1}{\partial s} & \dfrac{\partial h_2}{\partial s}\\ \dfrac{\partial h_1}{\partial t} & \dfrac{\partial h_2}{\partial t}\end{pmatrix}\right|_{(s,t,a)=(-(2a_0)/b^2\hat{t}_0, \hat{t}_0, a_0)}$$

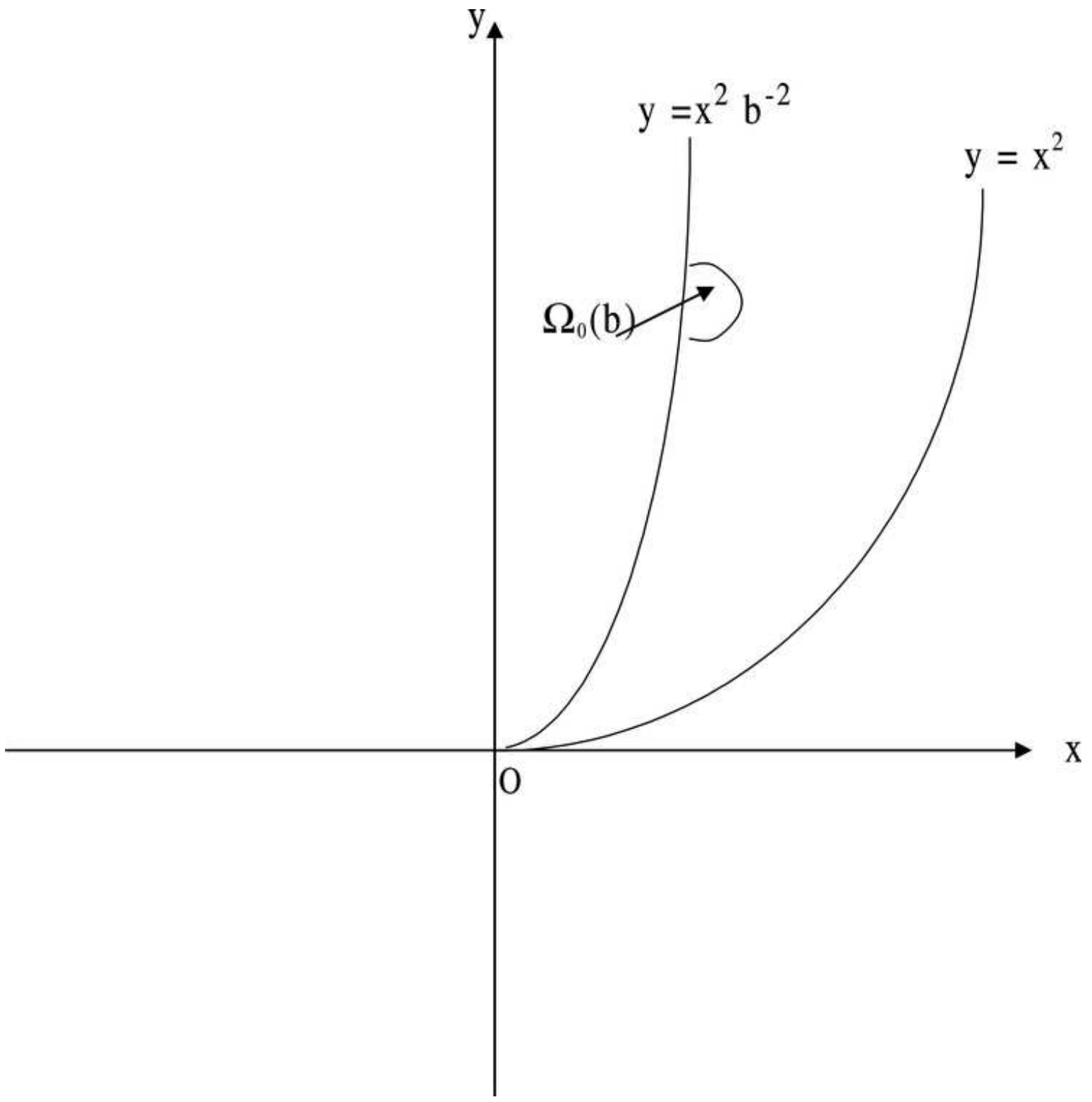

FIG. 5. *Partition of the area of integration.*

is positive definite, it follows from the implicit function theorem that there exists $\varepsilon > 0$ such that $\hat{s}_1 = \hat{s}(a, b_1)$ and $\hat{t}_1 = \hat{t}(a, b_1)$ are differentiable functions of $a$ and $b_1$ when $(a, a^2/b_1^2)^T \in \Omega_0(b)$ for any $0 < b < 1$, where $\hat{s}_1$ and $\hat{t}_1$ are solutions to the equations $K_s(s,t) = a, K_t(s,t) = a^2/b_1^2$. Since $\hat{t}_0 < 0$, we can always choose $\varepsilon$ to be so small that $\hat{t}_1 < 0$.

Using the transformation (3.5) and the saddlepoint approximation for $f_{(\overline{X},\overline{Y})}(x,y)$, we get

$$
\begin{aligned}
J_1(\tilde{b}) &= \iint_{\Omega_0(\tilde{b})} \hat{f}_n(x,y)\left(1 + \frac{r_n}{n}\right) dx\, dy \\
&= \iint_{(a,a^2/b^2)^T \in \Omega_0(\tilde{b})} \hat{f}_n(x(a,b), y(a,b))\, dx\, dy (1 + O(n^{-1})) \\
&= \iint_{(a,a^2/b^2)^T \in \Omega_0(\tilde{b})} \frac{n}{2\pi} \frac{\exp\{-n\Lambda(a,b)\}}{\det\{\Delta(a,b)\}^{1/2}} J(a,b)\, da\, db (1 + O(n^{-1})) \\
&= \int_{\tilde{b}}^{\tilde{b}+\delta_1} \int_{a_0-\delta_2}^{a_0+\delta_2} \frac{n}{2\pi} \frac{\exp\{-n\Lambda(a,b)\}}{\det\{\Delta(a,b)\}^{1/2}} J(a,b)\, da\, db (1 + O(n^{-1})),
\end{aligned} \tag{6.3}
$$

where $|r_n| < C$ since $\Omega_0(\tilde{b})$ is compact and $\delta_1$ and $\delta_2$ are small positive numbers such that

$$\text{if } a \in [a_0 - \delta_2, a_0 + \delta_2] \quad \text{and} \quad b \in [\tilde{b}, \tilde{b} + \delta_1] \qquad \text{then } (a, a^2/b^2)^T \in \Omega_0(\tilde{b}).$$

By Proposition 6.2, applying the Laplace approximation to the inner integral of (6.3) w.r.t. $a$ gives

$$\begin{aligned} J_1(\tilde{b}) &= \int_{\tilde{b}}^{\tilde{b}+\delta_1} \sqrt{\frac{n}{2\pi}} \frac{\exp\{-n\Lambda(a_0(b), b)\}}{\det\{\Delta(a_0(b), b)\}^{1/2}} \frac{J(a_0(b), b)}{\Lambda_{aa}^{1/2}(a_0(b), b)} \\ &\qquad \times \left(1 + \frac{r_{1n}}{n}\right) db(1 + O(n^{-1})) \qquad (6.4) \\ &= \int_{\tilde{b}}^{\tilde{b}+\delta_1} \sqrt{\frac{n}{2\pi}} \frac{\exp\{-n\Lambda(a_0(b), b)\}}{\det\{\Delta(a_0(b), b)\}^{1/2}} \frac{J(a_0(b), b)}{\Lambda_{aa}^{1/2}(a_0(b), b)} db(1 + O(n^{-1})), \end{aligned}$$

where $|r_{1n}|$ is uniformly bounded in $[\tilde{b}, \tilde{b} + \delta_1]$.

From Lemma A.8, $\Lambda(a_0(b), b)$ is a strictly increasing function of $b$ in the neighborhood of $b$. Define

$$\begin{aligned} w \equiv w(b) &= \sqrt{2\Lambda(a_0(b), b)}, \\ v \equiv v(b) &= \frac{\det\{\Delta(a_0, b)\}^{1/2} \Lambda_b(a_0, b) \Lambda_{aa}^{1/2}(a_0, b)}{|J(a_0, b)|}. \end{aligned}$$

Noting that

$$\frac{dw(b)}{db} = \frac{1}{w}\left(\Lambda_b(a_0, b) + \Lambda_a(a_0, b) \frac{da_0(b)}{db}\right) = \frac{\Lambda_b(a_0, b)}{w},$$

we have

$$J_1(\tilde{b}) = \int_{\tilde{w}}^{\tilde{w}_1} \sqrt{\frac{n}{2\pi}} \frac{e^{-nw^2/2}}{v} dw(1 + O(n^{-1})), \qquad (6.5)$$

where $\tilde{w} = w(\tilde{b})$ and $\tilde{w}_1 = w(\tilde{b} + \delta_1)$. Write $\tilde{v} = v(\tilde{b})$. Applying the Laplace approximation to the second integral of the following equality, we get

$$\begin{aligned} J_1(\tilde{b}) &= \int_{\tilde{w}}^{\tilde{w}_1} \sqrt{\frac{n}{2\pi}} e^{-nw^2/2} dw(1 + O(n^{-1})) \\ &\qquad - \int_{\tilde{w}}^{\tilde{w}_1} \sqrt{\frac{n}{2\pi}} e^{-nw^2/2} w\left(\frac{1}{w} - \frac{1}{v}\right) dw(1 + O(n^{-1})) \\ &= \Phi(\sqrt{n}\tilde{w}_1) - \Phi(\sqrt{n}\tilde{w}) - \frac{\phi(\sqrt{n}\tilde{w})}{\sqrt{n}}\left(\frac{1}{\tilde{w}} - \frac{1}{\tilde{v}} + O(n^{-1})\right) \qquad (6.6) \\ &= (1 - \Phi(\sqrt{n}\tilde{w}))(1 + O(n^{-1})) - \frac{\phi(\sqrt{n}\tilde{w})}{\sqrt{n}}\left(\frac{1}{\tilde{w}} - \frac{1}{\tilde{v}} + O(n^{-1})\right) \end{aligned}$$

$$= 1 - \Phi(\sqrt{n}\tilde{w}) - \frac{\phi(\sqrt{n}\tilde{w})}{\sqrt{n}}\left(\frac{1}{\tilde{w}} - \frac{1}{\tilde{v}} + O(n^{-1})\right),$$

where, in going from the second-to-the-last to the last line, we used $1 - \Phi(x) \sim \phi(x)/x$ as $x \to \infty$. Replacing $\tilde{b}$ by $b$, we get the desired result. □

PROPOSITION 6.6. *Under the conditions of Theorem 4.1,*

$$J_2(b)/J_1(b) = o(n^{-m}) \qquad \textit{for any } m > 0. \tag{6.7}$$

PROOF. By Lemma A.8, $\Lambda(a_0(b), b)$ is a strictly increasing function of $b$ for $b \in (0,1)$. Therefore, applying Laplace approximations to (6.4) again, we have

$$C_1 n \exp[-n\Lambda(a_0(b), b)] \le J_1(b) \le C_2 n \exp[-n\Lambda(a_0(b), b)]$$
$$\text{where } 0 < C_1 \le C_2 < \infty.$$

The proposition then follows from this and Lemma A.9. □

Finally, Theorem 4.1 follows from (6.1), (6.6) and (6.7).

## APPENDIX: SOME USEFUL LEMMAS

From here on, let $X$ be a r.v. with $EX = 0$ or $EX^2 = \infty$. We shall also adopt the same notation from Section 3. Write

$$I(s, t; a, b) = sa + ta^2/b^2 - K(s, t).$$

We now give our first lemma.

LEMMA A.1. *For fixed $a$ and $b$, we have*

$$\Lambda(a, b) = \sup_{s,t} I(s, t; a, b).$$

*When no solutions to $\partial I(s, t; a, b)/\partial s = \partial I(s, t; a, b)/\partial t = 0$ exist, we define $\Lambda(a, b) = \infty$.*

PROOF. It is easy to see that, for fixed $a$ and $b$, $I(s, t; a, b)$ is a concave function of $s$ and $t$ and it is differentiable for any $(s, t)^T \in \text{interior}(\Theta)$, where

$$\Theta = \{\theta = (s, t)^T : K(s, t) = \ln E e^{sX + tX^2} < \infty\}.$$

Therefore,

$$\sup_{s,t} I(s, t; a, b) = \hat{s}a + \hat{t}a^2/b^2 - K(\hat{s}, \hat{t}\,) = \Lambda(a, b),$$

where $\hat{s} = \hat{s}(x, y)$ and $\hat{t} = \hat{t}(x, y)$ are solutions to

$$K_s(\hat{s}, \hat{t}\,) = a, \qquad K_t(\hat{s}, \hat{t}\,) = a^2/b^2, \tag{A.1}$$

whenever the solutions exist. When no solutions exist, then clearly we have $\sup_{s,t} I(s,t;a,b) = \infty$. The proof is complete. $\square$

From Theorem 3.1 and Lemma A.1, we see that the saddlepoint approximation of the self-normalized sum involves finding, for fixed $b$,

$$\Lambda(a_0, b) := \inf_a \Lambda(a,b) = \inf_a \sup_{s,t} I(s,t;a,b) = I(\hat{s}_0, \hat{t}_0; a_0, b),$$

where $\hat{s}_0$, $\hat{t}_0$ and $a_0$ satisfy (3.6). In particular, we notice that the point $(\hat{s}_0, \hat{t}_0, a_0)^T$ falls on the curve $\hat{s}_0 = -2a_0\hat{t}_0/b^2$. This motivates the following definition:

$$g(t,a;b) = I(s,t;a,b)|_{s=-2at/b^2} = -ta^2/b^2 - K(-2at/b^2, t). \tag{A.2}$$

Also note that the domain of $a$ in the above infimum can be reduced to $\{a : ab > 0\}$ because of the transformation $a = x$ and $b = x/\sqrt{y}$. Since we only consider the case $0 < b < 1$, from now on we can suppose $a > 0$.

Let $C_s$ denote the support of the r.v. $X$, that is,

$$C_s = \{x : P(X \in (x - \varepsilon, x + \varepsilon)) > 0 \text{ for any } \varepsilon > 0\}.$$

Clearly, $C_s$ must be closed. We further use $\mathrm{Card}(C_s)$ to denote the number of elements in $C_s$ and define $\mathrm{Card}(C_s) = \infty$ if $C_s$ does not contain a finite number of elements.

LEMMA A.2. *Assume* $\mathrm{Card}(C_s) \geq 3$. *Then* $g(t,a;b)$ *is strictly decreasing in* $t$ *for* $t \in (-\varepsilon_0, \infty)$ *for some* $\varepsilon_0 > 0$.

PROOF. If suffices to show that $g(t,a;b)$ is strictly decreasing in $t$, either:

(I) for $t \in [0, \infty)$, or
(II) for $t \in (-\varepsilon_0, 0]$ for some $\varepsilon_0 > 0$.

We shall prove (I) first. Let $Z = -2aX/b^2 + X^2$. For arbitrary $t$ and $t_1$ such that $0 \leq t < t_1$, we need to show that $g(t,a;b) > g(t_1,a;b)$. If $Ee^{t_1 Z} = \infty$, then $g(t_1,a;b) = -a^2 t_1/b^2 - \ln Ee^{t_1 Z} = -\infty$, in which case (I) follows straightaway. Now, assume that $Ee^{t_1 Z} < \infty$ below, which implies that moments of $X$ of all orders exist. Thus, $g(t,a;b)$ is differentiable in $t$ for $t \in (-\infty, t_1)$. Taking derivatives gives

$$\frac{\partial g(t,a;b)}{\partial t} = -\frac{a^2}{b^2} - \frac{EZe^{tZ}}{Ee^{tZ}}. \tag{A.3}$$

Observe that

$$\left.\frac{\partial g(t,a;b)}{\partial t}\right|_{t=0} = -\frac{a^2}{b^2} - EX^2 < 0$$

and

$$\frac{\partial^2 g(t,a;b)}{\partial t^2} = -\left(\frac{EZ^2 e^{tZ}}{Ee^{tZ}} - \left(\frac{EZe^{tZ}}{Ee^{tZ}}\right)^2\right) < 0, \tag{A.4}$$

since $Z = X^2 - 2aX/b^2$ is nondegenerate by the assumption that $\mathrm{Card}(C_s) \geq 3$. Thus, $\frac{\partial g(t,a;b)}{\partial t} < 0$ when $t \in [0, t_1)$. So $g(t,a;b)$ is strictly decreasing in $[0, t_1)$. Since $t_1$ is arbitrary, we have hence proved (I).

We shall prove (II) next. If there exists some $t_2 > 0$ such that $Ee^{t_2 Z} < \infty$, then (II) follows from the fact that $\frac{\partial g(0,a;b)}{\partial t} = -a^2/b^2 - EX^2 < 0$. It remains to prove (II) under the condition that

$$Ee^{t_3 Z} = \infty \qquad \text{for all } t_3 > 0.$$

To show this, we choose an arbitrary $t < 0$. Then, from (A.3) we have

$$\frac{\partial g(t,a;b)}{\partial t} = -\frac{a^2}{b^2} - \frac{\int_{-\infty}^{\infty}(-2ax/b^2 + x^2)e^{t(x-a/b^2)^2}\,dF(x)}{\int_{-\infty}^{\infty} e^{t(x-a/b^2)^2}\,dF(x)}. \tag{A.5}$$

By the monotone convergence theorem we have

$$\begin{aligned} &\lim_{t\to 0^-}\int_{-\infty}^{\infty} e^{t(x-a/b^2)^2}\,dF(x) = 1, \\ &\lim_{t\to 0^-}\int_{-\infty}^{\infty} x^2 e^{t(x-a/b^2)^2}\,dF(x) = EX^2 \qquad \text{(maybe } \infty\text{)}, \end{aligned} \tag{A.6}$$

where $t \to 0^-$ means that $t \to 0$ from the left side of 0.

If $E|X| < \infty$, then noting $|xe^{t(x-a/b^2)^2}| \leq |x|$ for $t < 0$, we can use Lebesgue's dominated convergence theorem to get

$$\lim_{t\to 0^-}\int_{-\infty}^{\infty}\left(-\frac{2a}{b^2}x\right)e^{t(x-a/b^2)^2}\,dF(x) = -\frac{2a}{b^2}EX = 0. \tag{A.7}$$

If $E|X| = \infty$ (hence $EX^2 = \infty$), then

$$\begin{aligned} &\lim_{t\to 0^-}\int_{-\infty}^{\infty}\left(-\frac{2a}{b^2}x + x^2\right)e^{t(x-a/b^2)^2}\,dF(x) \\ &\quad \geq \lim_{t\to 0^-}\int_{-\infty}^{\infty}\left(-\frac{4a^2}{b^4} - \frac{x^2}{4} + x^2\right)e^{t(x-a/b^2)^2}\,dF(x) \\ &\quad = -\frac{4a^2}{b^4} + \frac{3}{4}EX^2 \\ &\quad = \infty. \end{aligned} \tag{A.8}$$

Combining (A.5)–(A.8) gives

$$\lim_{t\to 0^-}\frac{\partial g(t,a;b)}{\partial t} < 0.$$

Note that $g(t,a;b)$ is left continuous at $t = 0$. We conclude (II). $\square$

LEMMA A.3. *Assume that $F(x)$ is a continuous d.f. For each fixed $b \in (0,1)$ and $a \in R$, we have*

$$\sup_{t\in R} g(t,a;b) = \sup_{t<0} g(t,a;b), \tag{A.9}$$

*and the supremum is either attained at some finite unique point, $\tilde{t} := \tilde{t}(a,b) < 0$, or is simply infinity.*

PROOF. Define $h(x) := x^2 - 2ax/b^2 + a^2/b^2 = (x-a_1)(x-a_2)$, where

$$\begin{aligned} a_{10} &:= a_{10}(a) = \frac{a}{b^2}(1-\sqrt{1-b^2}\,), \\ a_{20} &:= a_{20}(a) = \frac{a}{b^2}(1+\sqrt{1-b^2}\,), \\ a_1 &:= a_1(a) = \min(a_{10}, a_{20}), \\ a_2 &:= a_2(a) = \max(a_{10}, a_{20}). \end{aligned} \tag{A.10}$$

Consider the following two cases:

(I′) $(a_1,a_2) \cap C_s \neq \varnothing$,
(II′) $(a_1,a_2) \cap C_s = \varnothing$.

First suppose that (I′) holds. Then there must exist $W := [a_3,a_4] \subset (a_1,a_2)$ so that:

(i) there exists $\delta > 0$ such that $h(x) < -\delta$ for each $x \in W$;
(ii) $P(X \in W) > 0$.

Then we have, as $t \to -\infty$,

$$\begin{aligned} g(t,a;b) &= -\ln \int_{-\infty}^{\infty} e^{th(x)}\,dF(x) \le -\ln \int_W e^{th(x)}\,dF(x) \\ &\le -\ln \int_W e^{-t\delta}\,dF(x) = t\delta - \ln P(X \in W) \\ &\to -\infty. \end{aligned}$$

From Lemma A.2, $\sup_{t\in R} g(t,a;b)$ is attained at some finite $\tilde{t} = \tilde{t}(a,b) < 0$. Since $g(t,a;b)$ is a differentiable function of $t$ when $t<0$, we have $\frac{\partial g(\tilde{t},a;b)}{\partial t} = 0$. This, together with (A.4), implies that there is at most one solution to the equation $\frac{\partial g(t,a;b)}{\partial t} = 0$. Therefore, $\tilde{t}$ is also unique.

Next suppose (II′) holds. Since $C_s$ is necessarily closed, then $[a_1,a_2] \cap C_s$ contains at most two points, $\{a_1,a_2\}$. Clearly, we have:

(i) $h(x) > 0$ for each $x \in C_s \setminus \{a_1,a_2\}$;
(ii) $P(X \in C_s \setminus \{a_1,a_2\}) > 0$,

where (ii) follows since $F(x)$ is a continuous d.f. Therefore, as $t \to -\infty$, we have

$$g(t,a;b) = -\ln \int_{-\infty}^{\infty} e^{th(x)}\, dF(x) = -\ln \int_{C_s \setminus \{a_1,a_2\}} e^{th(x)}\, dF(x) \to \infty. \qquad \square$$

Remark A.1. Lemma A.3 also holds true for $b \geq 1$, in which case both sides of (A.9) are equal to infinity.

Lemma A.4. *For $0 < b < 1$, define*

$$U = \{a : (a_1(a), a_2(a)) \cap C_s \neq \varnothing\},$$

*where $a_1(a)$ and $a_2(a)$ are defined in (A.10). Then, if $F(x)$ is a continuous d.f.:*

1. *$U$ is an open set and $U \neq \varnothing$, so does $U \cap R^+$, where $R^+ = \{x : x > 0\}$.*
2. *When $a \in U$, then $g(\tilde{t}(a,b), a; b) = \sup_{t<0} g(t,a;b) < \infty$, where $\tilde{t} = \tilde{t}(a,b) < 0$ is a finite unique solution to the equation $\frac{\partial g(t,a;b)}{\partial t} = 0$.*
3. *When $a \notin U$, then $\sup_{t<0} g(t,a;b) = \infty$.*
4. *$\inf_{a>0} \sup_{t \in R} g(t,a;b) = \inf_{a \in U \cap R^+} \sup_{t<0} g(t,a;b)$.*

Proof. We only prove 1 since 2–4 follow easily from the proof of Lemma A.3.

First, the claim that $U \neq \varnothing$ can be easily seen from the fact that $\bigcup_a \{a : (a_1(a), a_2(a))\} = R$. Second, we shall show that $U$ is open, which is equivalent to showing that the complement of $U$,

$$U_0 = \{a : (a_1(a), a_2(a)) \cap C_s = \varnothing\},$$

is a closed set. To show this, for any fixed $a' \in U_0$, then $(a_1(a'), a_2(a')) \not\subset C_s$, or $(a_1(a'), a_2(a')) \subset \overline{C}_s$, the complement of $C_s$. Let $V(a')$ be the largest interval such that $(a_1(a'), a_2(a')) \subset V(a') \subset \overline{C}_s$. For simplicity, assume that $a' > 0$ (the cases for $a' = 0$ and $a' < 0$ can be treated similarly). Since $\overline{C}_s$ is open, then $V(a')$ must be open as well. Write $V(a') = (c_0, d_0)$, where the endpoints could be $\infty$ or $-\infty$. Write

$$a_c(a') := \frac{c_0(1 + \sqrt{1-b^2}\,)}{b^2},$$

$$a_d(a') := \frac{d_0(1 - \sqrt{1-b^2}\,)}{b^2}.$$

It is easy to see that the closed interval $[a_c(a'), a_d(a')]$ will be the largest subset of $U_0$ including $a'$. Furthermore, for any $a' \neq a''$, the two intervals $[a_c(a'), a_d(a')]$ and $[a_c(a''), a_d(a'')]$ either coincide or are nonoverlapping. Therefore,

$$U_0 = \bigcup_{a' \in R} [a_c(a'), a_d(a')],$$

which is closed. The proof is complete. $\square$

LEMMA A.5. *For* $0 < b < 1$*:*

1. $\lim_{a\to\infty}\sup_{t<0} g(t,a;b) = \infty$, $\lim_{a\to 0^+}\sup_{t<0} g(t,a;b) = \infty$*;*
2. $\lim_{a\to\infty}\sup_{s\in R,t\in R} I(s,t;a,b) = \infty$, $\lim_{a\to 0^+}\sup_{s\in R,t\in R} I(s,t;a,b) = \infty$,

*where* $a \to 0^+$ *means that* $a$ *goes to* $0$ *from the right side.*

PROOF. Let $k$ be a positive number. Then

$$\begin{aligned}\sup_{t<0} g(t,a;b) &\geq g\left(-\frac{k}{a^2},a;b\right)\\ &= \frac{k}{b^2} - \ln\int_{-\infty}^{\infty}\exp\left\{-\frac{k}{a^2}\left(x^2-\frac{2a}{b^2}x\right)\right\}dF(x) \qquad \text{(A.11)}\\ &:= \frac{k}{b^2} - \ln M(a).\end{aligned}$$

It follows from Lebesgue's dominated convergence theorem that

$$\lim_{a\to\infty} M(a) = 1, \qquad \lim_{a\to 0^+} M(a) = 0. \qquad \text{(A.12)}$$

Combining (A.11) and (A.12) gives

$$\liminf_{a\to\infty}\sup_{t<0} g(t,a;b) \geq \frac{k}{b^2},$$

$$\liminf_{a\to 0^+}\sup_{t<0} g(t,a;b) = \infty.$$

Since $k$ can be arbitrarily large, we have proved 1.

From (A.2) we have $\sup_{s\in R,t\in R} I(s,t;a,b) \geq \sup_{t<0} g(t,a;b)$. This, together with 1 above, implies that

$$\lim_{a\to\infty}\sup_{s\in R,t\in R} I(s,t;a,b) = \infty,$$

$$\lim_{a\to 0^+}\sup_{s\in R,t\in R} I(s,t;a,b) = \infty,$$

which completes the proof of 2. □

LEMMA A.6. *Assume that* $F(x)$ *is a continuous d.f. and that* $0 < b < 1$*. Then* $\inf_{a>0}\sup_{t\in R} g(t,a;b)$ *is attained at some finite unique point,* $(a,t)^T = (a_0,\hat{t}_0)^T$*, where* $a_0 > 0$*,* $\hat{t}_0 := \tilde{t}(a_0,b) < 0$ *and they satisfy* (3.7).

PROOF. It follows from Lemmas A.3–A.5 that $\inf_{a>0}\sup_{t\in R} g(t,a;b)$ is attained at some finite points $a_0 \in U$ and $\hat{t}_0 := \tilde{t}(a_0,b) < 0$. When $a \in U$, by part 2 of Lemma A.4, we have

$$\frac{\partial g(\tilde{t},a;b)}{\partial t} = \frac{-EZe^{\tilde{t}Z}}{Ee^{\tilde{t}Z}} - \frac{a^2}{b^2} = 0, \qquad \text{where } Z = -\frac{2a}{b^2}X + X^2. \qquad \text{(A.13)}$$

By the assumption that $F(x)$ is a continuous d.f., which implies that $Z$ is nondegenerate, (A.4) is true. It then follows from the implicit function theorem that $\tilde{t}(a,b)$ is a differentiable function in some neighborhood $U^*(a)$ of $a$ (also a differentiable function in some neighborhood of $b$). We can also guarantee that $U^*(a) \subset U$. Hence $\sup_{t\in R} g(t,a;b)$ is also a differentiable function in some neighborhood of $a_0$. Thus $a_0$ satisfies the equation $\frac{dg(\tilde{t},a;b)}{da} = 0$, that is,

$$EX\exp\left\{\tilde{t}\left(-\frac{2a}{b^2}X + X^2\right)\right\} = aE\exp\left\{\tilde{t}\left(-\frac{2a}{b^2}X + X^2\right)\right\}. \tag{A.14}$$

It follows from (A.13) and (A.14) that $a_0$ and $\hat{t}_0$ are the solutions to the equations

$$\begin{aligned} EZe^{tZ} &= -\frac{a^2}{b^2}Ee^{tZ}, \\ EXe^{tZ} &= aEe^{tZ}, \end{aligned}$$

which are equivalent to (3.7) or (3.6).

Now we show the uniqueness of $(a_0,\hat{t}_0)^T$. Suppose $(a_0',\hat{t}_0')^T$ is another point such that $g(\hat{t}_0',a_0';b) = \inf_{a>0}\sup_{t\in R} g(t,a;b)$. Note that

$$g(t,a;b) = -\log E\exp\left\{t\left(-\frac{2a}{b^2}X + X^2 + \frac{a^2}{b^2}\right)\right\}.$$

We must have

$$\begin{aligned} E\exp\left\{\hat{t}_0\left(-\frac{2a_0}{b^2}X + X^2 + \frac{a_0^2}{b^2}\right)\right\} &= \sup_{a>0} E\exp\left\{\hat{t}_0\left(-\frac{2a}{b^2}X + X^2 + \frac{a^2}{b^2}\right)\right\} \\ &\geq E\exp\left\{\hat{t}_0\left(-\frac{2a_0'}{b^2}X + X^2 + \frac{a_0'^2}{b^2}\right)\right\} \\ &\geq \inf_{t<0} E\exp\left\{t\left(-\frac{2a_0'}{b^2}X + X^2 + \frac{a_0'^2}{b^2}\right)\right\} \\ &= E\exp\left\{\hat{t}_0'\left(-\frac{2a_0'}{b^2}X + X^2 + \frac{a_0'^2}{b^2}\right)\right\}. \end{aligned} \tag{A.15}$$

If $\hat{t}_0 \neq \hat{t}_0'$, then

$$E\exp\left\{\hat{t}_0\left(-\frac{2a_0'}{b^2}X + X^2 + \frac{a_0'^2}{b^2}\right)\right\} > \inf_{t<0} E\exp\left\{t\left(-\frac{2a_0'}{b^2}X + X^2 + \frac{a_0'^2}{b^2}\right)\right\} \tag{A.16}$$

by the fact that $E\exp\{t(-2aX/b^2 + X^2 + a^2/b^2)\}$ is a strictly convex function of $t$ for each fixed $a$ and $-\frac{2a}{b^2}X + X^2 + \frac{a^2}{b^2}$ is not identically equal to 0. Combining (A.15) and (A.16), we get

$$E\exp\left\{\hat{t}_0\left(-\frac{2a_0}{b^2}X + X^2 + \frac{a_0^2}{b^2}\right)\right\} > E\exp\left\{\hat{t}_0'\left(-\frac{2a_0'}{b^2}X + X^2 + \frac{a_0'^2}{b^2}\right)\right\},$$

which contradicts our assumption. Hence $\hat{t}_0 = \hat{t}'_0$.

Next we show that $\hat{s}_0 = \hat{s}'_0$. Define $f(a,s) = E\exp\{s(-2X/(ab^2) + X^2/a^2 + 1/b^2)\}$. Note that $f(a,s)$ is a strictly convex function of $s$ for each fixed $a$. Thus, we have

$$f(a_0, \hat{s}_0) = f(a'_0, \hat{s}'_0) = \sup_{a>0} \inf_{s<0} f(a,s),$$

where $\hat{s}_0 = \hat{t}_0 a_0^2$ and $\hat{s}'_0 = \hat{t}'_0 a_0'^2$. Similar to the proof of $\hat{t}_0 = \hat{t}'_0$ above, we can show that $\hat{s}_0 = \hat{s}'_0$. Hence $a_0 = a'_0$. This completes the proof of uniqueness. □

The next lemma establishes the relationship between $I(s,t;a,b)$ and $g(t,a;b)$.

Lemma A.7. *Assume that $F(x)$ is a continuous d.f. Then, for $0 < b < 1$,*

$$\inf_{a\ge 0} \sup_{t\le 0} g(t,a;b) = \inf_{a>0} \sup_{t<0} g(t,a;b) = \inf_{a>0} \sup_{s\in R, t\in R} I(s,t;a,b) \equiv \inf_{a>0} \Lambda(a,b).$$

Proof. The first equality holds since $g(t,a;b)$ is strictly decreasing as $t \to 0^-$ by Lemma A.2, and $\sup_{t<0} g(t,0;b) = \infty$. We shall now prove the second equality. From (A.2) we have

$$\inf_{a>0} \sup_{s\in R, t\in R} I(s,t;a,b) \ge \inf_{a>0} \sup_{t<0} g(t,a;b). \tag{A.17}$$

From 2 of Lemma A.5 we see that $\inf_{a>0} \sup_{s\in R, t\in R} I(s,t;a,b)$ is attained at some finite $\hat{a} > 0$. By Lemma A.6, $\inf_{a>0} \sup_{t<0} g(t,a;b)$ is also attained at some $a_0 > 0$ and $\hat{t}_0 < 0$ satisfying equation (3.7), namely,

$$K_s(-2a_0\hat{t}_0/b^2, \hat{t}_0) = a_0, \qquad K_t(-2a_0\hat{t}_0/b^2, \hat{t}_0) = a_0^2/b^2. \tag{A.18}$$

Therefore,

$$\begin{aligned}
&\inf_{a>0} \sup_{s\in R, t\in R} I(s,t;a,b) \\
&\quad = \sup_{s\in R, t\in R} I(s,t;\hat{a},b) \\
&\quad \le \sup_{s\in R, t\in R} I(s,t;a_0,b) \\
&\quad = \sup_{s\in R, t\in R} \{sa_0 + ta_0^2/b^2 - K(s,t)\} \\
&\quad = \{sa_0 + ta_0^2/b^2 - K(s,t)\}|_{s=\hat{s}_0, t=\hat{t}_0} \\
&\qquad\qquad [\text{where } K_s(\hat{s}_0,\hat{t}_0) = a_0,\ K_t(\hat{s}_0,\hat{t}_0) = a_0^2/b^2] \\
&\quad = \{sa_0 + ta_0^2/b^2 - K(s,t)\}|_{s=-2a_0\hat{t}_0/b^2, t=\hat{t}_0} \qquad [\text{by (A.18)}]
\end{aligned}$$

$$
\begin{aligned}
&= -\hat{t}_0 a_0^2/b^2 - K(-2a_0\hat{t}_0/b^2, \hat{t}_0) \\
&= g(\hat{t}_0, a_0; b) \\
&= \inf_{a>0} \sup_{t<0} g(t, a; b).
\end{aligned}
$$

The lemma thus follows from this and (A.17). □

LEMMA A.8. *Assume that $F(x)$ is a continuous d.f. Then, for $0<b<1$, $\inf_{a>0} \sup_{t<0} g(t,a;b)$ is a strictly increasing function of $b$.*

PROOF. Regard $g(\tilde{t}(a,b), a; b)$ as a joint function of $a$ and $b$. Then $\frac{\partial g(\tilde{t}(a,b),a;b)}{\partial b}|_{a=a_0} = -\frac{2a_0^2\tilde{t}(a_0,b)}{b^3} > 0$, that is, $g(\tilde{t}(a_0,b), a_0; b)$ is a strictly increasing function of $b$ in a small neighborhood of $b$. If $b_1 < b_2$ and $b_1$ is sufficiently close to $b_2$, we have

$$
g(\tilde{t}(a_1,b_1), a_1; b_1) \le g(\tilde{t}(a_2,b_1), a_2; b_1) < g(\tilde{t}(a_2,b_2), a_2; b_2), \tag{A.19}
$$

where $a_1$ and $a_2$ satisfy $g(\tilde{t}(a_1,b_1), a_1; b_1) = \inf_{a>0} g(\tilde{t}(a,b_1), a; b_1)$ and $g(\tilde{t}(a_2,b_2), a_2; b_2) = \inf_{a>0} g(\tilde{t}(a,b_2), a; b_2)$, respectively. Lemma A.7 and Proposition 6.3 imply that $\inf_{a>0} \sup_{t<0} g(t,a;b)$ is a nondecreasing function of $b$, which, combined with (A.19) holding under the condition that $b_1 < b_2$ and $b_1$ is sufficiently close to $b_2$, proves Lemma A.8. □

LEMMA A.9. *Assume that $F(x)$ is a continuous d.f. Then, for $0<b<1$, $\varepsilon > 0$ and $m > 0$, we have*

$$
P((\overline{X}, \overline{V}_n^2)^T \in (\Omega_1(b))^-)/\exp(-ng(\hat{t}_0, a_0; b)) = o(n^{-m}).
$$

PROOF. From Corollary 1.1 of Dembo and Shao (1998), we get

$$
\begin{aligned}
&\limsup_{n\to\infty} \frac{1}{n} \ln P((\overline{X}, \overline{V}_n^2)^T \in (\Omega_1(b))^-) \\
&\qquad \le - \inf_{(a, a^2/b_1^2)^T \in (\Omega_1(b))^-} \sup_{s,t} I(s,t;a,b_1) =: -I_{\min},
\end{aligned}
$$

if the condition (1.12) in Dembo and Shao (1998) holds, which is clearly true since

$$
\liminf_{y\to\infty, (x,y)^T \in (\Omega_1(b))^-} \frac{x^2}{y} = b^2 > 0.
$$

Hence, for all $\delta_1 > 0$, there exists $n_1$ such that if $n \ge n_1$,

$$
\frac{1}{n} \ln P((\overline{X}, \overline{V}_n^2)^T \in (\Omega_1(b))^-) \le -I_{\min} + \frac{\delta_1}{2}. \tag{A.20}
$$

From (A.2) we have $0 \le \sup_{t<0} g(t,a;b_1) \le \sup_{s,t} I(s,t;a,b_1)$ for any $b_1$, which implies that

$$-I_{\min} \le -\inf_{(a,a^2/b_1^2)^T \in (\Omega_1(b))^-} \sup_{t<0} g(t,a;b_1). \tag{A.21}$$

Define

$$\delta_1 = \inf_{(a,a^2/b_1^2)^T \in (\Omega_1(b))^-} \sup_{t<0} g(t,a;b_1) - g(\hat{t}_0, a_0; b).$$

We shall now show that $\delta_1 > 0$. Similar to Lemma A.5, we can show that if $b \le b' \le 1$, then

$$\lim_{a\to 0^+, b_1 \to b'} g(\tilde{t}, a; b_1) = \infty,$$

$$\lim_{a\to \infty, b_1 \to b'} g(\tilde{t}, a; b_1) = \infty.$$

Hence, $\inf_{(a,a^2/b_1^2)^T \in (\Omega_1(b))^-} g(\tilde{t}, a; b_1)$ is attained at some finite $a_E > 0$ and $b \le b_E \le 1$. By Lemma A.8 $g(\hat{t}_0, a_0; b_1) = \inf_{a>0} \sup_{t<0} g(t,a;b_1)$ is a strictly increasing function of $b_1$. If $b_E > b$, we have

$$\begin{aligned}
&\inf_{\{(a,b_1)^T : (a,a^2/b_1^2)^T \in (\Omega_1(b))^-\}} \sup_{t<0} g(t,a;b_1) \\
&\quad = \inf_{\{(a,b_1)^T : (a,a^2/b_1^2)^T \in (\Omega_1(b))^-\}} g(\tilde{t}(a,b_1), a; b_1) \\
&\quad = g(\tilde{t}(a_E, b_E), a_E; b_E) \\
&\quad \ge \inf_{a>0} g(\tilde{t}(a, b_E), a; b_E) \\
&\quad > \inf_{a>0} g(\tilde{t}(a, b), a; b) \\
&\quad = \inf_{a>0} \sup_{t<0} g(t,a;b).
\end{aligned}$$

By Lemma A.6, $a_0$ is unique. If $b_E = b$, we have

$$g(\tilde{t}(a_E, b_E), a_E; b_E) = g(\tilde{t}(a_E, b), a_E; b) > g(\hat{t}_0, a_0; b).$$

Combining the above facts, we have

$$\inf_{(a,a^2/b_1^2)^T \in (\Omega_1(b))^-} \sup_{t<0} g(t,a;b_1) > g(\hat{t}_0, a_0; b). \tag{A.22}$$

Therefore, we have proved that $\delta_1 > 0$. By (A.20)–(A.22), we have that if $n \ge n_1$,

$$P((\overline{X}, \overline{V}_n^2)^T \in (\Omega_1(b))^-) \le \exp\{-n g(\hat{t}_0, a_0; b) - n\delta_1/2\}.$$

The proof is complete. □

B.-Y. Jing
Department of Mathematics
Hong Kong University
of Science and Technology
Clear Water Bay, Kowloon
Hong Kong
e-mail: majing@ust.hk

Q.-M. Shao
Department of Mathematics
University of Oregon
Eugene, Oregon 97403
USA
and
Department of Mathematics
Department of Statistics
and Applied Probability
National University of Singapore
Singapore 117543
e-mail: qmshao@darkwing.uoregon.edu

W. Zhou
Department of Statistics
and Applied Probability
National University of Singapore
Singapore 117543
e-mail: stazw@nus.edu.sg